\input amstex
\documentstyle{amsppt}
\document
\magnification=1200
\NoBlackBoxes
\nologo
\vsize18cm

\centerline{\bf WEAK FROBENIUS MANIFOLDS}

\bigskip 

\centerline{C.~Hertling${}^1$, Yu.~Manin${}^2$}

\medskip

\centerline{{\it ${}^1$ Math.~Institut der Uni. Bonn, hertling\@math.uni-bonn.de}}

\centerline{{\it ${}^2$ MPIM, Bonn, Germany, manin\@mpim-bonn.mpg.de}}

\bigskip

{\bf Abstract.}
We establish a new universal relation between the Lie
bracket and $\circ$--multiplication of tangent fields
on any Frobenius (super)manifold. We use this identity
in order to introduce the notion of ``weak Frobenius manifold''
which does not involve  metric as part of structure.
As another application, we show that the powers of an
Euler field generate (a half of) the Virasoro algebra
on an arbitrary, not necessarily semi--simple,
Frobenius supermanifold.

\bigskip

{\bf 0. Introduction.} B.~Dubrovin introduced and
thoroughly studied in [D] the notion of Frobenius
manifold. By definition, it is a structure
$(M,g,\circ )$ where $M$ is a manifold,  $\circ$ is an 
associative, commutative and $\Cal{O}_M$--bilinear
multiplication on the tangent sheaf $\Cal{T}_M,$
and $g$ is a flat metric on $M,$ invariant with respect
to $\circ .$ The main axiom is the local existence of
a function $\Phi$ (Frobenius potential) such that
the structure constants of $\circ$ in the basis $\partial_a$ of flat
local fields are given by the tensor of third
derivatives $A_{ab}{}^c=\Phi_{ab}{}^c$ with one index raised with
the help of $g.$

\smallskip

We start with establishing a new universal identity (1)
between the $\circ$--multipli\-ca\-tion and the Lie bracket.
It  follows formally from the Poisson (or Leibniz) 
identity, but is strictly weaker, and
the algebra of tangent fields on a Frobenius manifold is never
Poisson. For further comments see section 5. We show
that this identity encodes an essential part of the
potentiality property, at least in the semisimple case.

\smallskip

We then use  it in order to introduce in section 3 ``weak Frobenius
manifolds'' that is, Frobenius manifolds without
a fixed flat metric. We explain the relation of 
this notion to  Dubrovin's notion of twisted
Frobenius manifolds ([D], Appendix B.) The importance
of weak Frobenius manifolds is related to the fact
that in the constructions of K.~Saito and Barannikov--Kontsevich
the metric is  the part of the structure that comes last,
and (at least in the theory of singularities) requires
considerable additional work.

\smallskip

Finally, in section 6 we extend the construction
of the Virasoro algebra from the Euler field, previously 
known only in the semisimple case,
to the general situation.

\smallskip

As a general reference on the basics of the theory
of Frobenius manifolds we use [M] (summarized in [MM].)
In particular, our manifolds are supermanifolds
(say, in the complex analytic category).
The multiplication $\circ$ is called semisimple,
if locally $(\Cal{T}_M,\circ )$ is isomorphic
to $\Cal{O}_M^n$ with componentwise multiplication.
The basic idempotent vector fields are then denoted
$e_i.$

\smallskip

We would like to thank Markus Rosellen for stimulating
discussions.

\medskip

\proclaim{\quad 1. Definition} An $F$--manifold
is a pair $(M,\circ ),$ where $M$
is a (super)manifold and $\circ$ is an associative
supercommutative $\Cal{O}_M$--bilinear multiplication
$\Cal{T}_M\times \Cal{T}_M\to \Cal{T}_M$ satisfying the
following identity: for any (local) vector fields
$X,Y,Z,W$ we have
$$
[X\circ Y,Z\circ W]-[X\circ Y,Z]\circ W-
(-1)^{(X+Y)Z}Z\circ [X\circ Y,W]
$$
$$
-X\circ [Y,Z\circ W] + X\circ [Y,Z]\circ W
+(-1)^{YZ} X\circ Z\circ [Y,W]
$$
$$
-(-1)^{XY}Y\circ [X,Z\circ W] + (-1)^{XY}Y\circ [X,Z]\circ W
+(-1)^{X(Y+Z)}Y\circ Z\circ [X,W] =0
\eqno(1)
$$
Here and in section 5 we write, say, $(-1)^{(X+Y)Z}$ as a shorthand for
$(-1)^{(\widetilde{X}+\widetilde{Y})\widetilde{Z}}$, where 
$\widetilde{X}$ is the parity of $X.$
\endproclaim

\medskip

{\bf 1.1. Remarks.} (i) The left hand side of (1)
is $\Cal{O}_M$--polylinear in $X,Y,Z,W$. In other words,
it is a tensor. This can be checked by a completely
straightforward, although lengthy, calculation.

\smallskip

(ii) Introduce the expression measuring the deviation
of the structure $(\Cal{T}_M, \circ , [,\,])$ from that
of Poisson algebra on  $(\Cal{T}_M, \circ )$:
$$
P_X(Z,W):= [X,Z\circ W]-[X,Z]\circ W-(-1)^{XZ}Z\circ [X,W].
\eqno(2)
$$
Then (1) is equivalent to the following requirement:
$$
P_{X\circ Y}(Z,W)=X\circ P_Y(Z,W)+(-1)^{XY}Y\circ P_X(Z,W).
\eqno(3)
$$

\medskip

\proclaim{\quad 2. Theorem} a). Let $(M,g,A)$  be a Frobenius
manifold with multiplication $\circ$. Then $(M,\circ )$
is an $F$--manifold.

\smallskip

b). Let $(M,\circ )$ be a pure even $F$--manifold,
whose multiplication law is semisimple on an open dense subset.
Assume that it admits an invariant flat metric $g$
defining  the cubic tensor $A$ as in [M], I.(1.2).
Then $(M,g,A)$ is Frobenius manifold.
\endproclaim

\medskip

{\bf Proof.} a). Since the left hand side of (1) is a tensor,
it suffices to check that it vanishes on quadruples of
flat fields $(X,Y,Z,W)=(\partial_a, \partial_b, \partial_c, \partial_d).$
Flat fields (super)commute so that only five summands of nine
survive in (1). Denoting the structure constants $A_{ab}{}^c$ as in [M],
I.(1.4) and calculating the coefficient of $\partial_f$ in (1),
we can represent it as a sum of five summands,
for which we introduce special notation in order
to explain the pattern of cancellation:
$$
\sum_e A_{ab}{}^e\partial_eA_{cd}{}^f-(-1)^{(a+b)(c+d)}
\sum_e A_{cd}{}^e\partial_eA_{ab}{}^f=\alpha_1+\beta_1,
$$
$$
(-1)^{(a+b)c}\sum_e \partial_cA_{ab}{}^eA_{ed}{}^f=\alpha_2,\quad
(-1)^{(a+b+c)d}\sum_e \partial_dA_{ab}{}^eA_{ec}{}^f=\gamma_1,
$$
$$
-(-1)^{a(b+c+d)}\sum_e \partial_bA_{cd}{}^eA_{ea}{}^f=\gamma_2,\quad
-(-1)^{(c+d)b}\sum_e \partial_aA_{cd}{}^eA_{eb}{}^f=\beta_2.
$$
Here we write, say, $(-1)^{(a+b)c}$ as a shorthand for
$(-1)^{(\widetilde{x_a}+\widetilde{x_b})\widetilde{x_c}}$.

\smallskip

Use potentiality in order to interchange the subscripts $e,c$
in $\alpha_1$ (a sign emerges). After this we see that
$$
\alpha_1+\alpha_2=(-1)^{(a+b)c}\partial_c
\left(\sum_e A_{ab}{}^eA_{ed}{}^f\right).
$$
Similarly, permuting $a$ and $e$ in $\beta_1$ we find
$$
\beta_1+\beta_2=-(-1)^{(c+d)b}\partial_a
\left(\sum_e A_{cd}{}^eA_{eb}{}^f\right).
$$ 
Now rewrite $\gamma_1$ permuting $a,d,$ and $\gamma_2$
permuting $b,c.$ Calculating finally 
$\beta_1+\beta_2+\gamma_1+\gamma_2$ we see that it cancels with
$\alpha_1+\alpha_2$ due to the associativity relations [M], I.(1.5).

\smallskip

b). Clearly, $(M,g,A)$ is an associative pre--Frobenius manifold
in the sense of [M], Definition I.1.3, so that it only remains
to check its potentiality in the domain of semisimplicity.
To this end we will use the Theorem I.3.3 of [M]. 

\smallskip

Let $(e_i)$ be the basic idempotent local
vector fields. Applying (3) to $X=Y=e_i$ we get
$P_{e_i}=2e_i\circ P_{e_i}$ so that $P_{e_i}=0.$
Applying then (2) to $(X,Z,W)=(e_i,e_j,e_j),\,i\ne j,$
we see that $[e_i,e_j]=0.$ This is the first condition
of the Theorem I.3.3. The second one expresses invariance and flatness of the
metric in canonical coordinates, which we have already
postulated. 

\medskip

\proclaim{\quad 2.1. Corollary (of the proof)} Semisimple
$F$--manifolds are exactly those manifolds $(M,\circ )$
which everywhere locally admit a basis of pairwise commuting
$\circ$--idempotent vector fields, or, which is
the same, Dubrovin's canonical coordinates.
\endproclaim 

\smallskip

In fact, we have already deduced from
(1) that that  $e_i$ pairwise commute.
Conversely, if they commute, (1) holds for any quadruple
of idempotents.

\medskip

\proclaim{\quad 3. Definition} A weak Frobenius manifold
is an $F$--manifold $(M,\circ )$ such that in a neighborhood $U$
of any point there exists a flat invariant metric $g$
making $(U,g,\circ )$ Frobenius manifold. We will
call such metrics compatible (with the given $F$--structure.)
\endproclaim

\smallskip

Thus, a weak Frobenius manifold is a Frobenius
manifold without a fixed metric.

\smallskip

Any semisimple $F$--manifold is automatically
weak Frobenius. This follows from the results
of [M], Chapter II, \S 3, which reduce the construction
of compatible metrics to the solution to
Schlesinger's equations. We do not know whether
there exist non--semisimple $F$--manifolds which are not 
weak Frobenius.

\medskip

{\bf 4. Sheaves of compatible metrics and Euler fields.} 
Let $(M,\circ )$ 
be a weak Frobenius manifold.
Then compatible metrics on $M$ form a sheaf $\Cal{M}_M$. 
Assume now that $M$ admits an identity $e.$
Then there is an embedding $\Cal{M}_M\hookrightarrow \Omega^1_M$
which sends each metric $g$ to the respective coidentity
$\varepsilon_g$ defined by $i_X\varepsilon_g=g(e,X).$ In fact,
knowing $\varepsilon_g$ we can reconstruct 
$g:\, g(X,Y)=i_{X\circ Y}(\varepsilon_g).$ We will call
$\varepsilon_g$ compatible 1--forms.

\smallskip

We will now impose an additional restriction and
denote by $\Cal{F}_M$ the sheaf of those compatible metrics
for which $e$ is flat. We will call such metrics {\it admissible.}
If $g\in \Cal{F}_M,$ then
$\varepsilon_g$ is closed and $g$--flat: see [M], I.(2.4).
It is important to understand the structure of
this sheaf of sets. Again, the situation is rather
transparent on the tame semisimple part of $M.$ 
We will state and prove
Dubrovin's theorem which provides a neat local
description of admissible metrics {\it with fixed
rotation coefficients $\gamma_{ij}$} considered as
functions on the common definition domain of metrics.

\smallskip

This result should be compared with the Theorem II.3.4.3
of [M] which depicts  the set of metrics {\it with fixed
$v_{ij}$ at a point} where
$$
v_{ij}=\frac{1}{2}\,(u^j-u^i)\,\frac{\eta_{ij}}{\eta_j},\ 
\gamma_{ij} = \frac{1}{2}\,\frac{\eta_{ij}}{\sqrt{\eta_i\eta_j}}=
\frac{1}{u^j-u^i}\sqrt{\frac{\eta_j}{\eta_i}}\,v_{ij}
\eqno(4)
$$
so that both statements refer to the closely related coordinates
on the space of metrics.

\medskip

\proclaim{\quad 4.1. Theorem} a). Let $g =\sum_i\eta_i (du^i)^2,\,
\widetilde{g} =\sum_i\widetilde{\eta}_i (du^i)^2$ be two 
$\circ$--invariant metrics
in a simply connected domain of canonical coordinates $u^1,\dots ,u^n$ in $M.$
Then there exist exactly $2^n$ vector fields $\partial$ in this
domain such that 
$$
\widetilde{g}(X,Y)=g(\partial\circ X,\partial\circ Y)
\eqno(5)
$$
for any $X,Y.$ These fields are $\circ$--invertible
and differ only by the signs of their $e_i$--components.

\smallskip

b). Assume that $g$ is admissible. Then $\widetilde{g}$
defined by (5) is admissible and has the same
rotation coefficients $\widetilde{\gamma}_{ij}=\gamma_{ij}$ 
iff $\partial$ is
$g$--flat.
\endproclaim

\smallskip

{\bf Proof.} a). Put $\partial =\sum_iD_ie_i.$ Then
(5) is equivalent to 
$$
D_i^2=\frac{\widetilde{\eta}_i}{\eta_i}.
\eqno(6)
$$
This proves the first statement.

\smallskip

b). Choose a solution $(D_i)$ of (6). Let $\nabla_i$ denote
the Levi--Civita covariant derivative in the direction $e_i$
with respect to the metric $g.$ Using [M], I.(3.10),
we see that
$$
\nabla_i(\sum_jD_je_j)=
$$
$$
\sum_{j\ne i}\left( e_iD_j+\frac{1}{2}\,D_j\frac{\eta_{ij}}{\eta_j}-
\frac{1}{2}\,D_i\frac{\eta_{ij}}{\eta_j}\right) e_j
+\left(e_iD_i+\frac{1}{2}\sum_j D_j\frac{\eta_{ij}}{\eta_i}\right) e_i.
\eqno(7)
$$
Using (6) and [M], I.(3.13), we can rewrite the first sum in (7) as
$$
\sum_{j\ne i}\left( \sqrt{\frac{\widetilde{\eta}_i}{\eta_j}}
(\widetilde{\gamma}_{ij}-\gamma_{ij})\right)\,e_j \,
\eqno(8)
$$
and the remaining terms as
$$
\left(\frac{1}{2}\,\frac{\widetilde{\eta}_{ii}}{\sqrt{\eta_i\widetilde{\eta}_i}}
+\sum_{j\ne i}\sqrt{\frac{\widetilde{\eta_j}}{\eta_i}}\,\gamma_{ij}
\right) e_i\,.
\eqno(9)
$$
If we replace in (9) $\gamma_{ij}$ by $\widetilde{\gamma}_{ij}$,
the resulting expression will vanish for admissible
$\widetilde{g}$ because $(\sum_je_j)\widetilde{\eta}_i=0$
(see [M], Proposition I.3.3). Subtracting this zero from (9)
we finally find
$$
\nabla_i(\sum_jD_je_j)=
\sum_{j\ne i} \sqrt{\frac{\widetilde{\eta}_i}{\eta_j}}
(\widetilde{\gamma}_{ij}-\gamma_{ij})\, e_j
-\sum_{j\ne i} \sqrt{\frac{\widetilde{\eta}_j}{\eta_i}}
(\widetilde{\gamma}_{ij}-\gamma_{ij})\, e_i \,.
\eqno(10)
$$
Hence admissibility
of $\widetilde{g}$ and coincidence of the rotation coefficients 
imply the $g$--flatness of $\partial ,$ and vice versa. 
This proves the second statement of the theorem.

\smallskip

If one does not assume semi--simplicity, a part of
the preceding theorem still holds true. It is, too, due
to Dubrovin.

\medskip

\proclaim{\quad 4.2. Theorem} Let $g$ be an admissible metric and
$\partial$ a $g$--flat even invertible vector field.
Then $\widetilde{g}$ defined by (5) is admissible.
\endproclaim

\smallskip

{\bf Proof.} Put $\widetilde{x}^a:=\sum_b g^{ab}\partial_b\partial\Phi 
=\partial\Phi^a,$
where $(\partial_a=\partial /\partial x^a)$ is a local basis of flat vector 
fields and $\Phi$ is a local Frobenius potential.
Then we have
$$
\frac{\partial \widetilde{x}^a}{\partial x^b}=
 \partial \Phi_b^a.
$$
Hence the respective Jacobian is the matrix of the
multiplication $\partial\circ$ in the basis $(\partial_a).$
Since the latter is invertible, $(\widetilde{x}^a)$ 
form a local coordinate system. For the dual
basis of vector fields $\widetilde{\partial}_a$
we have
$$
\partial\circ\widetilde{\partial}_a=
\sum_b \widetilde{\partial}_a(x^b)\partial_b\circ\partial=
\sum_{b,c,d} \widetilde{\partial}_a(x^b)g^{cd}\partial_b\partial_d
\partial\Phi\partial_c=
\sum_{b,c} \widetilde{\partial}_a(x^b)\partial_b(\widetilde{x}^c)\partial_c=
\partial_a\, .
\eqno(11)
$$
Hence $\widetilde{g}$ has the same coefficients in the
basis $(\widetilde{\partial}_a)$ as $g$ in the basis 
$(\partial_a)$ and is flat. Clearly, it is also
$\circ$--invariant. It remains to show that
the pre--Frobenius structure $(M,\circ ,\widetilde{g})$
is potential. It is easy to see that any local function
$\widetilde{\Phi}$ satisfying the equations
$$
\widetilde{\partial}_a\widetilde{\partial}_b\widetilde{\Phi}=
\partial_a\partial_b\Phi
\eqno(12)
$$
for all $a,b$ can serve as a local potential defining
the same multiplication $\circ$. To prove its
existence, we check
the integrability condition:
$$
\widetilde{\partial}_a\partial_b\partial_c\Phi=
\sum_d \widetilde{\partial}_a(x^d)\partial_d\partial_b\partial_c\Phi=
\sum_d \widetilde{\partial}_a(x^d)g(\partial_d,\partial_b\circ\partial_c)
$$
$$
=g(\widetilde{\partial}_a,\partial_b\circ\partial_c)=
g(\partial^{-1}\circ\partial_a,\partial_b\circ\partial_c)=
(-1)^{ab}\widetilde{\partial}_b\partial_a\partial_c\Phi\,.
\eqno(13)
$$
The same reasoning as in the end of the proof of [M], Theorem I.1.5
then shows the existence of $\widetilde{\Phi}.$ The identity
$e$ remains flat because $\widetilde{g}$--flat fields
form the sheaf $\partial^{-1}\Cal{T}^f_M$ and $e=\partial^{-1}\circ\partial .$
This finishes the proof.

\smallskip

Dubrovin calls the passage from $g$ to $\tilde{g}$
{\it the Legendre--type transformation.} In the Appendix B
of [D] he also constructs a different type of transformations
which he calls {\it inversion.}

\smallskip

What Dubrovin calls {\it a twisted Frobenius manifold}
in our language is a weak Frobenius manifold,
endowed with local admissible metrics connected
by the Legen\-dre--type transformations on the
overlaps of their definition domains.

\smallskip

We now turn to Euler fields. Let again $(M,\circ )$ be
a weak Frobenius manifold. An even vector field $E$ on M is called {\it a
weak Euler
field of (constant) weight $d_0$} if $\roman{Lie}_E(\circ )=d_0\circ$ that is,
for all local vector fields $X,Y$ we have
$$
P_E(X,Y)=[E,X\circ Y]-[E,X]\circ Y-X\circ [E,Y]=d_0X\circ Y.
\eqno(14)
$$
This is the same as [M], I.(2.6). If $(M,\circ )$ admits an identity $e,$
we get formally from (14) that $[e,E]=d_0e.$
Clearly, local weak Euler fields form a sheaf of vector spaces $\Cal{E}_X$,
and weight is a linear function on this sheaf. 
If $(M,\circ )$ comes from a Frobenius manifold with flat
identity $e$, then any Euler vector field
on the latter is a weak Euler field, and $e$ itself is a (weak)
Euler field of weight zero. The latter statement
follows by combining [M], Proposition I.2.2.2 and I.(2.3).

\medskip

\proclaim{\quad 4.3. Proposition} Commutator of
any two (local) weak Euler fields is a weak Euler field of
weight zero.
\endproclaim

\smallskip

{\bf Proof.} We start with the following general identity: for any
local vector fields $X,Y,Z,W$ we have
$$
P_{[X,Y]}(Z,W)=
$$
$$
[X,P_Y(Z,W)]-(-1)^{XY}P_Y([X,Z],W)-(-1)^{X(Y+Z)}P_Y(Z,[X,W])
$$
$$
-(-1)^{XY}[Y,P_X(Z,W)]+P_X([Y,Z],W)+(-1)^{YZ}P_X(Z,[Y,W])\,.
\eqno(15)
$$
In order to check this, replace the seven terms in (15)
by their expressions from (2), and then rewrite
the resulting three terms in the left hand side
using the Jacobi identity. All the twenty four summands will cancel.

\smallskip

Now apply (15) to the two weak Euler fields $X=E_1,\,Y=E_2.$
The right hand side will turn to zero. This proves our statement.

\medskip

{\bf 4.4. Example (Sh. Katz).} The (formal) Frobenius
manifold corresponding to the quantum cohomology
of a projective algebraic manifold $V$ admits at least
two different Euler fields (besides $e$), if $h^{pq}(V)\ne 0$
for some $p\ne q.$ To write them down explicitly,
choose a basis $(\partial_a )$ of $H=H^*(V,\bold{C})$
considered as the space of flat vector fields, and let
$(x^a)$ be the dual flat coordinates vanishing at zero. 
Let $\partial_a\in H^{p_a,q_a}(V).$ 
Put $-K_V=\sum_{p_b+q_b=2}r^b\partial_b.$ Then 
$$
E_1:=\sum_a (1-p_a)x^a\partial_a+\sum_br^b\partial_b,
$$
$$
E_2:=\sum_a (1-q_a)x^a\partial_a+\sum_br^b\partial_b
$$
are Euler.

\medskip

Let now $g$ be an admissible metric on $(M,\circ ).$
Weak Euler fields which are conformal with respect
to $g$ form a subsheaf of linear spaces in $\Cal{E}_M$
endowed with a linear function $D$, conformal weight:
see [M], I.(2.5). A direct calculation shows that
the commutator of such fields is conformal
of conformal weight zero. One can also say what happens to
the weights (and the full spectrum) of $E$ when
one replaces $g$ by another metric as in (5).

\medskip
 
\proclaim{\quad 4.5. Proposition} Let $(M,\circ ,g)$ be
a Frobenius manifold with flat identity $e$ and an Euler field $E$,
$[e,E]=d_0e, \roman{Lie}_E(g)=Dg, (d_a) =$ the spectrum of
$-\roman{ad}\,E$ on flat vector fields.
Assume that $\partial$ is an invertible flat field
such that $[\partial ,E]=d\partial.$ Then $E$ is an Euler field on 
$(M,\widetilde{g},\circ ),\roman{Lie}_E(\widetilde{g})=
(D+2d_0-2d)\widetilde{g},$ and the spectrum of $-\roman{ad}\,E$
on $\partial^{-1}\Cal{T}^f_M$ is $(d_a+d_0-d).$  
\endproclaim

\smallskip

We leave the straightforward 
proof to the reader.

\medskip

{\bf  5. Relation to Poisson structure.} Consider an abstract
structure $(A,\circ , [,])$ where $\circ$, resp. $[,]$
induce on the $\bold{Z}_2$--graded additive group $A$
the structure of supercommutative, resp. Lie, ring.
Assume that these operations satisfy the relation
(1), or equivalently, (3). Then we will call
$(A,\circ ,[,]),$ or simply $A$, {\it an $F$--algebra.} In particular,
vector fields on a Frobenius manifold form a sheaf
of $F$--algebras. 

\smallskip

Every Poisson algebra is an $F$--algebra. Conversely,
let $A$ be an $F$--algebra, and
$$
B:=\{X\in A\,|\,P_X\equiv 0\}.
\eqno(16)
$$

\smallskip

\proclaim{\quad 5.1. Proposition} a). $B$ is closed with respect
to $\circ$ and $[,]$ and hence forms a Poisson subalgebra.
If $A$ contains identity $e,$ then $e\in B.$

\smallskip

b). If $A$ is the algebra of vector fields on a split semisimple
Frobenius manifold, $B$ is spanned by the basic idempotent fields $e_i$
over constants. In particular, the Lie bracket in $B$ is trivial.
\endproclaim

\smallskip

{\bf Proof.} a). Assume that $P_X=P_Y=0$. We have 
$P_{X\circ Y}=0$ in view of
(3). Putting $X=Y=e$ in (3), we get $P_e=0.$
Finally,  $P_{[X,Y]}=0$ follows from (15).

\smallskip

b). Writing $X,Y,Z$ in the basis $e_i$ with indeterminate 
coefficients, one easily checks that if $P_X(Y,Z)=0$ for all $Y,Z,$
then the coefficients of $X$ are constant.

\medskip

\proclaim{\quad 6. Theorem} Let $E$ be an Euler field
on a Frobenius manifold with identity $e$ such that
$[e,E]=d_0e.$ Then for all $m,n\ge 0$
$$
[E^{\circ n},E^{\circ m}]=d_0(m-n)\,E^{\circ m+n-1}.
\eqno(17)
$$
\endproclaim
\smallskip

{\bf Proof.} We will prove slightly more. Let $X$ be
an even vector field on an arbitrary $F$--manifold with identity $e.$
Since in view of (3) the map $X \mapsto P_X$ 
is a $\circ$--derivation, we have $P_{X^{\circ n}}=nX^{\circ n-1}\circ P_X.$
Moreover, from (2) we have
$$
P_{X^{\circ n}}(e,e)=-[X^{\circ n},e].
$$
Hence
$$
[X^{\circ n},e]=nX^{\circ n-1}\circ [X,e].
\eqno(18)
$$
Let us {\it assume} now that $X$ satisfies the following
identities: for all $n\ge 1$
$$
[X^{\circ n},X]=(1-n)\,X^{\circ n}\circ [e,X].
\eqno(19)
$$
Then we assert that for all $m,n\ge 0$
$$
[X^{\circ n},X^{\circ m}]=(m-n)\,X^{\circ n+m-1}\circ [e,X].
\eqno(20)
$$
In fact, the cases when $m$ or $n$ is $\le 1$
are covered by (18), (19). The general case
can be treated by induction. We have
$$
[X^{\circ n},X^{\circ m}]=P_{X^{\circ n}}(X^{\circ m-1},X)+
[X^{\circ n},X^{\circ m-1}]\circ X+
[X^{\circ n},X]\circ X^{\circ m-1}=
$$
$$
nX^{\circ n-1}\circ ([X,X^{\circ m}]-[X,X^{\circ m-1}]\circ X)+
[X^{\circ n},X^{\circ m-1}]\circ X+
[X^{\circ n},X]\circ X^{\circ m-1}
$$
$$
=(m-n)\,X^{\circ n+m-1}\circ [e,X].
$$
It remains to notice that since $[e,E]=d_0e,$
$E$ satisfies (19) in view  of the general identity [M], I.(2.12).

\medskip

{\bf 6.1. Remark.} In the semisimple case
the meaning of (19) is transparent:
writing $X=\sum X_ie_i,$ we must have $e_iX_j=0$ for $i\ne j.$

\bigskip

\centerline{\bf References}

\medskip

[D] B.~Dubrovin. {\it Geometry of 2D topological field theories.}
In: Springer LNM, 1620 (1996), 120--348

\smallskip

[M] Yu.~Manin {\it Frobenius manifolds, quantum cohomology,
and moduli spaces (Chapters I,\,II,\,III).} Preprint MPI  96--113,
1996.

\smallskip

[MM] Yu.~Manin, S.~Merkulov. {\it Semisimple Frobenius
(super)manifolds and quantum cohomology of $\bold{P}^r$.}
Topological Methods in Nonlinear Analysis, 9:1 (1997), 107--161,
alg--geom/9702014.

\enddocument